\documentclass[a4paper,conference]{IEEEtran}
\IEEEoverridecommandlockouts
\usepackage{cite}
\usepackage{amsmath,amssymb,amsfonts}
\usepackage{algorithmic}
\usepackage{graphicx}
\usepackage{textcomp}
\usepackage{xcolor}
\def\BibTeX{{\rm B\kern-.05em{\sc i\kern-.025em b}\kern-.08em
    T\kern-.1667em\lower.7ex\hbox{E}\kern-.125emX}}

\usepackage{amsmath,amsfonts,amssymb,amsthm,mathtools}
\usepackage{hhline}
\usepackage{clrscode}

\newtheorem{exmp}{Example}

\newcommand{\E}[1]{\mathop{{\rm \bf E}\!\left\{#1\right\}}\nolimits}

\title{One-Step Condensed Forms for Square-Root Maximum Correntropy Criterion Kalman Filtering
\thanks{
The author acknowledges the financial support of the Portuguese FCT~--- \emph{Funda\c{c}\~ao para a Ci\^encia e a Tecnologia}, through the project UID/Multi/04621/2019 of CEMAT/IST-ID, Center for Computational and Stochastic Mathematics, Instituto Superior T\'ecnico, University of Lisbon.}
}

\author{\IEEEauthorblockN{Maria V.~Kulikova}
\IEEEauthorblockA{\textit{Center for Computational and Stochastic Mathematics} \\
\textit{Instituto Superior T\'{e}cnico, Universidade de Lisboa}\\
Lisbon, Portugal \\
maria.kulikova at ist.utl.pt}
}

\begin{document}

\maketitle

\begin{abstract}
This paper suggests a few novel Cholesky-based square-root algorithms for the maximum correntropy criterion Kalman filtering. In contrast to the previously obtained results, new algorithms are developed in the so-called {\it condensed} form that corresponds to the {\it a priori} filtering. Square-root filter implementations are known to possess a better conditioning and improved numerical robustness when solving ill-conditioned estimation problems. Additionally, the new algorithms permit easier propagation of the state estimate and do not require a back-substitution for computing the estimate. Performance of novel filtering methods is examined by using a fourth order benchmark navigation system example.
\end{abstract}

\section{INTRODUCTION}\label{sect1}

The Kalman filtering (KF) like estimators developed under the maximum correntropy criterion (MCC) methodology are shown to outperform the classical KF for estimation quality in case of non-Gaussian and impulsive noise scenario in state-space models, both for linear and nonlinear systems~\cite{2007:Liu,2014:Chen,2015:Chen,2011:Cinar,2012:Cinar,2016:Izanloo,2017:Chen,2017:Liu:KF,2017:CSL:Kulikova,2017:Yang}. Among nonlinear MCC filtering strategies, we may mention the MCC extended KF in~\cite{2016:Liu:EKF,2016:Wang}, the MCC unscented KF in~\cite{2017:Liu:UKF,2017:Wang}, the MCC-based Gauss-Hermite quadrature filter in~\cite{2017:Qin} and the so-called accurate continuous-discrete extended KF~\cite{2016:Kulikov:SISCI,2017:Kulikov:ANM} explored under the MCC in~\cite{2018:Kulikov:SP}. Here, we focus on the MCC-KF estimators suggested in~\cite{2016:Izanloo,2017:CSL:Kulikova}.

The numerical stability issues of the cited MCC-KF estimators are of special interest in this paper. Recently, their reliable square-root implementation methods have been designed in~\cite{2019:SP:Kulikova}. Following the approaches existed in the KF community, the new MCC-KF implementations belong to three main {\it factored-form} classes: (i) the Cholesky-based algorithms; e.g.,~\cite{1963:Potter,1995:Park,2018:Kulikov:IET}, (ii) the $UDU^{\top}$ factorization-based methods, e.g.,~\cite{1976:Thornton:PhD,1977:Bierman}, and (iii) singular value decomposition solution, e.g.,~\cite{2017:Kulikova:IET,2018:Kulikova:Romania}. All mentioned square-root strategies have been utilized for the robust MCC-KF filtering in~\cite{2019:SP:Kulikova}. In the cited paper, the implementation methods are developed in terms of covariance quantities and in the {\it a posteriori} form that yields the two-stage implementations. In contrast to the previously obtained results, the goal of this paper is to suggest square-root MCC-KF algorithms in the {\it a priori} form that results to the {\it condensed} implementations. We stress that the stated problem has never been solved in engineering literature, before.

More precisely, we develop the {\it a priori} covariance filtering approach for the original MCC-KF proposed in~\cite{2016:Izanloo} and for its improved variant (IMCC-KF) published in~\cite{2017:CSL:Kulikova}. Additionally, we discuss the square-root problem for the examined estimators and suggest a few condensed implementations for the robust filtering in terms of {\it lower} triangular Cholesky factors. As all square-root implementations, the new algorithms possess a better conditioning, improved numerical robustness (with respect to roundoff errors) and reliability when solving ill-conditioned estimation problems. Additionally, the new algorithms permit easier propagation of the state estimate and do not require a back-substitution for computing the estimate. It is achieved by utilizing orthogonal rotations for propagating the involved Cholesky factors as far as possible. Finally, we stress that derivation of new implementation methods is important for its own sake because it provides practitioners a diversity of methods with a fair possibility to choose any of them depending on practical application, complexity and accuracy requirements. Performance of novel filtering algorithms is examined by using a
 fourth order benchmark navigation system example.

The paper is organized as follows. Section~\ref{problem:statement} contains the problem statement and presents the derivation of the {\it a priori} MCC-KF recursions. The new square-root condensed implementation methods can be found in Section~\ref{main:result} together with their derivation. Section~\ref{numerical:experiments} contains the results of numerical simulation study and Section~\ref{conclusion} concludes the paper.

\section{PROBLEM STATEMENT} \label{problem:statement}

Consider a linear discrete-time stochastic system
  \begin{align}
   x_{k+1} = & F_k x_{k}  + G_k w_{k}, \label{eq:st:1} \\
    y_k  = & H_k x_k + v_k , \quad k \ge 0  \label{eq:st:2}
  \end{align}
where system matrices $F_k \in \mathbb R^{n\times n}$, $G_k \in \mathbb R^{n\times q}$ and $H_k \in \mathbb R^{m\times n}$ are known and constant over time. The vectors $x_k \in \mathbb R^n$  and $y_k \in \mathbb R^m$ are unknown dynamic state and available measurement vector, respectively. We assume that random variables $x_0$, $w_k$ and $v_k$ satisfy
\begin{align*}
&\E{x_0}       = \bar x_0,    &  \!\!\!&\E{(x_0-\bar x_0)(x_0-\bar x_0)^{\top}}  = \Pi_0, \\
&\E{w_k}       = \E{v_k} = 0, &  \!\!\!&\E{w_kx_0^{\top}} = \E{v_kx_0^{\top}}  = 0,  \\
&\E{w_kv_k^{\top}}  = 0,           &  \!\!\!&\E{w_kw_j^{\top}} = Q_k\delta_{kj}, \\
&&&  \E{v_kv_j^{\top}}  = R_k\delta_{kj}
\end{align*}
where covariance matrices $Q_k \in \mathbb R^{q\times q}$ and $R_k \in \mathbb R^{m\times m}$ are known. The symbol $\delta_{kj}$ is the Kronecker delta function.

The classical KF produces the minimum {\it linear} expected mean square error (MSE) estimate of the unknown state vector $x_k$.
It can be formulated in the {\it a priori} form as shown in~\cite[Theorem~9.2.1]{2000:Kailath:book}; see also the derivation in~\cite[Section~5.3]{2006:Simon:book}. More precisely, the one-step ahead predicted estimate $\hat x_{k+1|k}$ ({\it a priori} estimate) is calculated by
\begin{align}
\hat x_{k+1|k}  & =   F_k \hat x_{k|k-1}+K_{p,k}(y_k-H_k\hat x_{k|k-1}),  \label{KF:X:Riccati} \\
 R_{e,k} & = R_k+H_kP_{k|k-1}H_k^{\top},     \label{KF:Rek:Riccati} \\
K_{p,k}  & = F_kP_{k|k-1}H_k^{\top}R_{e,k}^{-1}    \label{KF:K:Riccati}
\end{align}
where $e_k = y_k-H_k\hat x_{k|k-1}$ are innovations (residuals) of discrete-time KF with covariance $R_{e,k} = \E{e_ke_k^{\top}}$ and $K_{p,k}=\E{\hat x_{k+1|k} e_k^{\top}} = F_kP_{k|k-1}H_k^{\top}$. The matrix $P_{k|k-1}$ is one-step ahead predicted error covariance $P_{k|k-1}=\E{ (x_{k}-\hat
x_{k|k-1})(x_{k}-\hat x_{k|k-1})^{\top}}$ computed via the Riccati difference recursion:
\begin{equation}
P_{k+1|k}  = F_kP_{k|k-1}F_k^{\top} + G_kQ_kG_k^{\top} - K_{p,k}R_{e,k}K_{p,k}^{\top} \label{KF:P:Riccati}
\end{equation}
with initial values $P_{0|-1} =  \Pi_0 > 0$ and $\hat x_{0|-1}  = \bar x_0$.

Being a {\it linear} estimator, the classical KF exhibits only sub-optimal behavior in non-Gaussian settings. To enhance its estimation quality and robustness with respect to outliers (impulsive noise), the maximum correntropy criterion has been utilized together with the KF estimation cost function in~\cite{2012:Cinar,2016:Izanloo}. In general, the {\it correntropy} is a similarity measure between two random variables~\cite{2007:Liu}. Following~\cite[Chapter~5]{2018:Principe:book}, it can be used in the estimation problem as follows: an estimator of unknown state $X \in {\mathbb R}$ is defined as a function of measurements $Y \in {\mathbb R}^m$, i.e., $\hat X = g(Y)$ where $g$ is solved by maximizing the correntropy between $X$ and $\hat X$ that is~\cite{2012:Chen}
\begin{equation} \label{mcc:kriterion}
\mbox{arg}\max \limits_{g \in G} V(X,\hat X) = \mbox{arg}\max \limits_{g \in G} \E{k_{\sigma}\Bigl(X - g(Y)\Bigr)}
\end{equation}
where $G$ stands for the collection of all measurable functions of $Y$, $k_{\sigma}(\cdot)$ is a kernel function and  $\sigma > 0$ is the kernel size. The Gaussian kernel is the most widely used kernel function and it is given as follows:
\begin{equation}\label{Gauss_kernel}
k_{\sigma}(X - \hat X) = \exp \left\{ -{(X - \hat X)^2}/{(2\sigma^2)}\right\}.
\end{equation}
It is not difficult to see that the MCC cost~\eqref{mcc:kriterion} with Gaussian kernel~\eqref{Gauss_kernel} reaches its maximum if and only if $X = \hat X$.

In~\cite{2016:Izanloo}, the MCC-KF is developed by solving the following estimation problem with the Gaussian kernel:
\begin{align*}
\mbox{arg}\max & \left\{  k_{\sigma}(\|\hat x_{k|k}-F_{k-1}\hat x_{k-1|k-1}\|)  \right. \\
               &  \left. + k_{\sigma}(\|y_k-H_k\hat x_{k|k}\|) \right\}
\end{align*}
where weighted norms are additionally utilized.

Next, the resulted nonlinear equation arisen in the optimization problem above is solved by a fixed point rule that yields the following recursion~\cite[p.~503]{2016:Izanloo}:
\begin{equation}
\hat x_{k|k}  =  \hat x_{k|k-1}  + K_k^{\lambda}(y_k - H_k\hat x_{k|k-1}) \label{eq:result}
\end{equation}
where $\hat x_{k|k-1}:=F_{k-1}\hat x_{k-1|k-1}$ and the gain matrix is defined as $K_k^{\lambda} = \lambda_k (P_{k|k-1}^{-1}+\lambda_k H_k^{\top} R_k^{-1} H_k)^{-1}H_k^{\top} R_k^{-1}$ with a scalar adjusting weight $\lambda_k$ given as follows:
\begin{equation}
\lambda_{k}  = \frac{k_{\sigma}(\|y_k-H_k\hat x_{k|k-1}\|_{R_k^{-1}})}{k_{\sigma}(\|\hat x_{k|k-1}-F_{k-1}\hat x_{k-1|k-1}\|_{P_{k|k-1}^{-1}})}.  \label{eq:lambda}
\end{equation}
Finally, the resulted MCC-KF estimator in~\cite{2016:Izanloo} utilizes the recursion for state estimate in~\eqref{eq:result} together with a symmetric Joseph stabilized equation for calculating the error covariance matrix $P_{k|k}$ that is known to improve the numerical robustness of KF-like implementations because of ensured symmetric form of $P_{k|k}$ in the presence of roundoff errors~\cite{2006:Simon:book,2015:Grewal:book}. In summary, the MCC-KF estimator is given as follows~\cite{2016:Izanloo}.

{\bf MCC-KF estimator}. \textsc{Time Update} ($k=\overline{1,K}$).
\begin{align}
\hat x_{k|k-1} & = F_{k-1}\hat x_{k-1|k-1}, \label{mcckf:p:X} \\
P_{k|k-1} & = F_{k-1}P_{k-1|k-1}F_{k-1}^T+G_{k-1}Q_{k-1}G_{k-1}^{\top}. \label{mcckf:p:P}
\end{align}

\textsc{Measurement Update}. Find {\it a posteriori} estimate:
\begin{align}
K_{k}^{\lambda} & = \lambda_k \left( P_{k|k-1}^{-1}+\lambda_k H_k^{\top} R_k^{-1} H_k \right)^{-1}H_k^{\top} R_k^{-1}, \label{mcckf:f:K:old}  \\
\hat x_{k|k}    & = \hat x_{k|k-1}+K_{k}^{\lambda}e_k, \quad  e_k  = y_k-H_k\hat x_{k|k-1},   \label{mcckf:f:X} \\
P_{k|k}  & = (I- K_{k}^{\lambda}H_k)P_{k|k-1}(I- K_{k}^{\lambda}H_k)^{\top} \nonumber \\
& +K_k^{\lambda} R_k [K_k^{\lambda}]^{\top} \label{mcckf:f:P:old}
\end{align}
where the scaling (inflation) parameter $\lambda_{k}$ is computed by~\eqref{eq:lambda}.

In~\cite{2017:CSL:Kulikova}, the estimation quality of original MCC-KF estimator presented above has been boosted by deriving mathematically equivalent formulas for the gain and error covariance, similar to the KF equations in~\cite[pp.~128-129]{2006:Simon:book}. It yields the improved MCC-KF (IMCC-KF) estimator summarized below.

{\bf IMCC-KF estimator}. \textsc{Time Update} ($k=\overline{1,K}$). Calculate {\it a priori} estimate $\hat x_{k|k-1}$ and $P_{k|k-1}$ by equations~\eqref{mcckf:p:X}, \eqref{mcckf:p:P}.

\textsc{Measurement Update}. Find {\it a posteriori} estimate $\hat x_{k|k}$ by formula~\eqref{mcckf:f:X} with the scaling parameter $\lambda_{k}$ from~\eqref{eq:lambda} and
\begin{align}
K_{k}^{\lambda} & = \lambda_k P_{k|k-1}H_k^{\top}[R_{e,k}^{\lambda}]^{-1}, \label{mcckf:f:K}  \\
R_{e,k}^{\lambda} & = \lambda_kH_kP_{k|k-1}H_k^{\top}+R_k, \label{mcckf:f:Rek} \\
P_{k|k}  & = (I - K^{\lambda}_{k}H_k)P_{k|k-1}. \label{mcckf:f:P}
\end{align}

As can be seen, the original MCC-KF and IMCC-KF equations recursively updates  {\it a posteriori} estimate $\hat x_{k|k}$ given the initial value $\hat x_{0|0}:= \bar x_0$. The discussed estimators can be re-formulated in the so-called {\it a priori} form similar to the classical KF equations~\eqref{KF:X:Riccati}~-- \eqref{KF:P:Riccati}. This yields one-step KF-like equations and allows for convenient {\it condensed} form of the resulted square-root implementation methods. Next, we derive one-step equations for both the MCC-KF and IMCC-KF estimators and explore the possibility of designing their Cholesky factorization-based implementations.

\subsection{The MCC-KF estimator in a priori form}

Following~\cite[Section~5.3]{2006:Simon:book}, we derive {\it a priori} MCC-KF recursion related to the original MCC-KF equations~\eqref{mcckf:p:X}--\eqref{mcckf:f:P:old}. Having substituted~\eqref{mcckf:f:X} into~\eqref{mcckf:p:X}, we get the recursion for one-step ahead predicted ({\it a priori}) estimate $\hat x_{k+1|k}$ as follows:
\begin{align}
\hat x_{k+1|k}      & =   F_k\hat x_{k|k} = F_k(\hat x_{k|k-1}+K_{k}^{\lambda}e_k),  \nonumber \\
& = F_k\hat x_{k|k-1} + F_kK_k^{\lambda}e_k = F_k\hat x_{k|k-1} + K_{p,k}^{\lambda}e_k  \label{MCCKF:X:Riccati}
\end{align}
where we introduced the notation $K_{p,k}^{\lambda} = F_kK_k^{\lambda}$.

It is proved in~\cite{2017:CSL:Kulikova} that equations~\eqref{mcckf:f:K:old} and~\eqref{mcckf:f:K} for gain calculation $K_k^{\lambda}$ are equivalent. Hence, we have
\begin{align}
K_{p,k}^{\lambda} &  = \lambda_k F_k \left( P_{k|k-1}^{-1}+\lambda_k H_k^{\top} R_k^{-1} H_k \right)^{-1}H_k^{\top} R_k^{-1}\\
& = \lambda_kF_kP_{k|k-1}H_k^{\top}[R_{e,k}^{\lambda}]^{-1}.  \label{MCCKF:K:Riccati}
\end{align}

Having substituted formula~\eqref{mcckf:f:P:old} into~\eqref{mcckf:p:P} and taking into account formula~\eqref{MCCKF:K:Riccati} and $K_{p,k}^{\lambda} = F_kK_k^{\lambda}$, we obtain
\begin{align*}
P_{k+1|k} & = F_kP_{k|k}F_k^{\top}+G_kQ_kG_k^{\top} \\
&  = F_k(I- K_{k}^{\lambda}H_k)P_{k|k-1}(I- K_{k}^{\lambda}H_k)^{\top}F_k^{\top} \\
& +F_kK_k^{\lambda} R_k [K_k^{\lambda}]^{\top}F_k^{\top} + G_kQ_kG_k^{\top} \\
&  = (F_k- K_{p,k}^{\lambda}H_k)P_{k|k-1}(F_k- K_{p,k}^{\lambda}H_k)^{\top} \\
& + K_{p,k}^{\lambda} R_k [K_{p,k}^{\lambda}]^{\top} + G_kQ_kG_k^{\top} \\
&  = F_k P_{k|k-1}F_k^{\top} +  G_kQ_kG_k^{\top} - F_k P_{k|k-1}H_k^{\top}[K_{p,k}^{\lambda}]^{\top}\\
& -K_{p,k}^{\lambda}H_kP_{k|k-1}F^{\top} \\
& + K_{p,k}^{\lambda}(H_kP_{k|k-1}H_k^{\top}+R_k)[K_{p,k}^{\lambda}]^{\top}.
\end{align*}
Taking into account that $H_kP_{k|k-1}H_k^{\top}+R_k = R_{e,k}$, formula~\eqref{MCCKF:K:Riccati} and symmetric form of $R_{e,k}$ and $R_{e,k}^{\lambda}$, we get
\[F_kP_{k|k-1}H_k^{\top} = \lambda^{-1}_k K_{p,k}^{\lambda}R_{e,k}^{\lambda},\]
we obtain the recursion for {\it a priori} error covariance $P_{k+1|k}$:
\begin{align*}
P_{k+1|k} & = F_k P_{k|k-1}F_k^{\top} +  G_kQ_kG_k^{\top} - 2\lambda^{-1}_kK_{p,k}^{\lambda}R_{e,k}^{\lambda}[K_{p,k}^{\lambda}]^{\top}\\
& + K_{p,k}^{\lambda}R_{e,k}[K_{p,k}^{\lambda}]^{\top} = F_k P_{k|k-1}F_k^{\top} +  G_kQ_kG_k^{\top} \\
& - K_{p,k}^{\lambda}[2\lambda^{-1}_kR_{e,k}^{\lambda}-R_{e,k}][K_{p,k}^{\lambda}]^{\top} \\
& = F_k P_{k|k-1}F_k^{\top} +  G_kQ_kG_k^{\top} \\
& - K_{p,k}^{\lambda}\Bigl(H_kP_{k|k-1}H_k^{\top}+(2\lambda_k^{-1}-1)R_k\Bigr)[K_{p,k}^{\lambda}]^{\top}.
\end{align*}

As can be seen, the MCC-KF equation for matrix $P_{k+1|k}$ computation does not allow the related square-root form because the term $2\lambda_k^{-1}-1$ might be negative, although the adjusting weight $\lambda_k$ defined in~\eqref{eq:lambda} is a nonnegative value. Meanwhile, for the improved MCC-KF counterpart the efficient square-root implementations do exist as shown below.

\subsection{The IMCC-KF estimator in a priori form}

Similarly, {\it a priori} recursion for the IMCC-KF equations~\eqref{mcckf:p:X}, \eqref{mcckf:p:P}, \eqref{mcckf:f:X}, \eqref{mcckf:f:K}--\eqref{mcckf:f:P} has been recently derived in~\cite{2020:IEEE:Kulikova}. For readers' convenience, we briefly discuss it here. First, one-step ahead predicted estimate $\hat x_{k+1|k}$ is calculated by equation~\eqref{MCCKF:X:Riccati} as in the MCC-KF examined above.

Next, having substituted formula~\eqref{mcckf:f:P} into~\eqref{mcckf:p:P} and taking into account formula~\eqref{MCCKF:K:Riccati} and symmetric form of any covariance matrix, we obtain
\begin{align*}
P_{k+1|k} & = F_kP_{k|k}F_k^{\top}+G_kQ_kG_k^{\top} \\
&  = F_k\left[(I - K^{\lambda}_{k}H_k)P_{k|k-1}\right]F_k^{\top}+G_kQ_kG_k^{\top} \\
& = F_kP_{k|k-1}F_k^{\top}+G_kQ_kG_k^{\top} - F_kK^{\lambda}_{k}H_kP_{k|k-1}F_k^{\top} \\
& = F_kP_{k|k-1}F_k^{\top} +G_kQ_kG_k^{\top} - K^{\lambda}_{p,k}(F_kP_{k|k-1}H_k^{\top})^{\top} \\
& = F_kP_{k|k-1}F_k^{\top} +G_kQ_kG_k^{\top} - \frac{1}{\lambda_k}K^{\lambda}_{p,k} R_{e,k}^{\lambda}[K^{\lambda}_{p,k}]^{\top}
\end{align*}
or, alternatively, we can avoid a scalar division by introducing
\begin{align}
K_{p,k}  & = \frac{1}{\lambda_k} F_kK_k^{\lambda} = F_kP_{k|k-1}H_k^{\top}[R_{e,k}^{\lambda}]^{-1},  \label{MCCKF:Knew:Riccati}
\end{align}
i.e., $K_{p,k}^{\lambda} = \lambda_k K_{p,k}$ defined in~\eqref{MCCKF:K:Riccati}. Thus, the equation for error covariance computation is re-formulated as follows:
\begin{align}
\!\!P_{k+1|k}  = & F_kP_{k|k-1}F_k^{\top} +G_kQ_kG_k^{\top} - \lambda_kK_{p,k}R_{e,k}^{\lambda}K_{p,k}^{\top} \label{MCCKF:P:Riccati1}  \\
\mbox{where } & K_{p,k} = F_kP_{k|k-1}H_k^{\top}[R_{e,k}^{\lambda}]^{-1}, \label{MCCKF:K:Riccati1} \\
& \hat x_{k+1|k}  = F_k\hat x_{k|k-1} + \lambda_kK_{p,k}e_k. \label{MCCKF:X:Riccati1}
\end{align}

Formulas~\eqref{MCCKF:P:Riccati1}--\eqref{MCCKF:X:Riccati1} are the {\it a priori} form of the IMCC-KF estimator, which we are looking for. It is important that the adjusting weight $\lambda_k$ defined in~\eqref{eq:lambda} is a nonnegative value and, hence, a square root exists (for real nonnegative numbers).

\section{MAIN RESULTS: CONDENSED FORMS} \label{main:result}

The resulted one-step IMCC-KF equations~\eqref{MCCKF:P:Riccati1}~--~\eqref{MCCKF:X:Riccati1} allow for convenient {\it condensed} form of stable square-root implementation methods that we derive in this section. Our square-root solution is based on the Cholesky decomposition of a symmetric positive definite matrix $A$ in the form $A=A^{1/2}A^{{\top}/2}$ where the factor $A^{1/2}$ is a {\it lower} triangular matrix with positive diagonal elements.

The condensed form implies that  one-step IMCC-KF equations~\eqref{MCCKF:P:Riccati1}~--~\eqref{MCCKF:X:Riccati1} are summarized into unique pre-array and, next, a stable orthogonal rotation is applied for obtaining the resulted post-array with the propagated filter quantities. Such implementation methods are additionally suitable for parallel implementation. Thus, we formulate the first condensed implementation method for the IMCC-KF estimator.

\begin{codebox}
\Procname{{\bf Algorithm 1}. $\proc{SR IMCC-KF}$ ({\it a priori, one-step form})}
\zi \textsc{Initialization:}($k=0$)
\li \>Apply Cholesky factorization: $\Pi_0 = \Pi_0^{1/2}\Pi_0^{{\top}/2}$;
\li \>Set initial values: $\hat x_{0|-1} = \bar x_0$, $P_{0|-1}^{1/2} = \Pi_0^{1/2}$;
\zi \textsc{Filter Recursion}: ($k=\overline{0,K}$)
\li \>Compute $\lambda_k$ and residual $e_k$ by formulas~\eqref{eq:lambda}, \eqref{mcckf:f:X};
\li \>Build pre-array and block-triangularize it   \label{mcc:sr1:p:P}
\zi \>$\!\!\!\!\underbrace{
\left[
\begin{smallmatrix}
R_k^{1/2} & \lambda_k^{1/2}H_{k}P_{k|k-1}^{1/2} & 0 \\
0 & F_kP_{k|k-1}^{1/2} & G_{k}Q_{k}^{1/2}
\end{smallmatrix}
\right]
}_{\mbox{\scriptsize Pre-array ${\mathbb A}_k$}} {\mathbb Q}_k \!\!=
\underbrace{
\left[
\begin{smallmatrix}
[R_{e,k}^{\lambda}]^{1/2} & 0 & 0 \\
\bar K_{p,k}^{\lambda} & P_{k+1|k}^{1/2} & 0
\end{smallmatrix}
\right]
}_{\mbox{\scriptsize Post-array ${\mathbb R}_k$}}$;
\li \>Read-off from post-array: $P_{k+1|k}^{1/2}$, $[R_{e,k}^{\lambda}]^{1/2}$, $\bar K_{p,k}^{\lambda}$;
\li \>Find $\hat x_{k+1|k}  = F_k\hat x_{k|k-1} + \lambda_k^{1/2}[\bar K_{p,k}^{\lambda}][R_{e,k}^{\lambda}]^{-1/2} e_k$.   \label{mcc:sr1:p:X}
\end{codebox}

Following~\cite{1995:Park}, it is not difficult to prove that the transformation in line~\ref{mcc:sr1:p:P} of Algorithm~1 implies the calculation by formulas~\eqref{MCCKF:P:Riccati1} and~\eqref{MCCKF:K:Riccati1} of the conventional one-step implementation. Let's consider the pre-array in line~\ref{mcc:sr1:p:P} of Algorithm~1
\[
\begin{bmatrix}
R_k^{1/2} & \lambda_k^{1/2}H_{k}P_{k|k-1}^{1/2} & 0 \\
0 & F_kP_{k|k-1}^{1/2} & G_{k}Q_{k}^{1/2}
\end{bmatrix}
\]
and block-triangularize it by any orthogonal rotation, i.e.,
\[{\mathbb A}_k {\mathbb Q}_k =
\begin{bmatrix}
X & 0 & 0 \\
Y & Z & 0
\end{bmatrix}.
\]

It is worth noting here that the orthogonal transformation ${\mathbb Q}_k$ sets up a conformal (i.e., a norm- and angle-preserving) mapping between the (block) rows of the pre-array ${\mathbb A}_k$ and post-array. Thus, we have
\begin{align*}
 <[R_k^{1/2}, \lambda_k^{1/2}H_{k}P_{k|k-1}^{1/2}, 0], & [R_k^{1/2}, \lambda_k^{1/2}H_{k}P_{k|k-1}^{1/2}, 0]> \\
 = & <[X, 0, 0],[X, 0, 0]>, \\
 <[R_k^{1/2}, \lambda_k^{1/2}H_{k}P_{k|k-1}^{1/2}, 0], & [0, F_kP_{k|k-1}^{1/2}, G_{k}Q_{k}^{1/2}]> \\
 = & <[X, 0, 0],[Y, Z, 0]>, \\
 <[0, F_kP_{k|k-1}^{1/2}, G_{k}Q_{k}^{1/2}],& [0, F_kP_{k|k-1}^{1/2}, G_{k}Q_{k}^{1/2}]> \\
 = & <[Y, Z, 0],[Y, Z, 0]>.
\end{align*}

From the first equality above, we obtain
\[XX^{\top} = \lambda_kH_k(P_{k|k-1}^{1/2}P_{k|k-1}^{{\top}/2})H_k^{\top} + R_k.\]
Having compared with equation~\eqref{mcckf:f:Rek}, we conclude that  $XX^{\top} = R_{e,k}^{\lambda}$ and, hence, $X:= [R_{e,k}^{\lambda}]^{1/2}$.

At the same way, from the second equality we get
\[XY^{\top} = \lambda_k^{1/2}H_{k}(P_{k|k-1}^{1/2}P_{k|k-1}^{{\top}/2})F_k^{\top}\]
and, hence, we define $Y = \lambda_k^{1/2}F_kP_{k|k-1}H_{k}^T[R_{e,k}^{\lambda}]^{-{\top}/2}:=\bar K_{p,k}^{\lambda}$. Having compared the normalized quantity $\bar K_{p,k}^{\lambda}$ naturally appeared in square-root Algorithm~1 with the value $K_{p,k}$ defined by equation~\eqref{MCCKF:K:Riccati1} in the conventional one-step method, we obtain the relationship
\[\bar K_{p,k}^{\lambda} = \lambda_k^{1/2}K_{p,k}[R_{e,k}^{\lambda}]^{1/2} \mbox{ or } K_{p,k} = \lambda_k^{-1/2}\bar K_{p,k}^{\lambda}[R_{e,k}^{\lambda}]^{-1/2} \]

Taking into account the formula above, we express equation~\eqref{MCCKF:X:Riccati1} of the conventional method in terms of values available in the square-root Algorithm~1, i.e., in terms of $[\bar K_{p,k}^{\lambda}]$ and $[R_{e,k}^{\lambda}]^{1/2}$ that are directly read-off from the post-array
\begin{align*}
\hat x_{k+1|k} & = F_k\hat x_{k|k-1} + \lambda_kK_{p,k}e_k \\
& = F_k\hat x_{k|k-1} + \lambda_k^{1/2}[\bar K_{p,k}^{\lambda}][R_{e,k}^{\lambda}]^{-1/2} e_k
\end{align*}
that is exactly the formula in line~\ref{mcc:sr1:p:X} of Algorithm~1.

Finally, we consider the last equality, i.e.,
\[F_k(P_{k|k-1}^{1/2}P_{k|k-1}^{{\top}/2})F^{\top}_k+G_kQ_kG_k^{\top} =  YY^{\top}+ZZ^{\top}.
\]
Having substituted $Y=\bar K_{p,k}^{\lambda}$ and $\bar K_{p,k}^{\lambda} = \lambda_k^{1/2}K_{p,k}[R_{e,k}^{\lambda}]^{1/2}$, we get
\begin{align*}
ZZ^{\top} & = F_kP_{k|k-1}F_k^{\top} +G_kQ_kG_k^{\top} - \bar K_{p,k}^{\lambda} [\bar K_{p,k}^{\lambda}]^{\top} \\
& = F_kP_{k|k-1}F_k^{\top} +G_kQ_kG_k^{\top} - \lambda_kK_{p,k}R_{e,k}^{\lambda}K_{p,k}^{\top}
\end{align*}
i.e., from~\eqref{MCCKF:P:Riccati1} we conclude $Z:=P_{k+1|k}^{1/2}$. Thus, the square-root Algorithm~1 is justified.

The newly-derived square-root condensed form in Algorithm~1 is equivalent to conventional one-step formulas~\eqref{MCCKF:P:Riccati1}~--~\eqref{MCCKF:X:Riccati1}, but it is convenient for practical application and more suited to parallel computations. Besides, the square-root formulation makes the new Algorithm~1 is inherently more stable (with respect to roundoff errors) than any conventional one-step implementation; e.g., the one-step formulas~\eqref{MCCKF:P:Riccati1}~--~\eqref{MCCKF:X:Riccati1}.

We further suggest one more stable square-root implementation developed in the condensed form with the following additional benefit. As can be seen,
Algorithms~1 allows for simple propagation of the filter quantities $[\bar K_{p,k}^{\lambda}]$ and $[R_{e,k}^{\lambda}]^{1/2}$ at each iteration step. Indeed, these values are simply read-off from the resulted post-array obtained from block-triangularization of the pre-array. Although stable orthogonal rotations are used for $[\bar K_{p,k}^{\lambda}]$ and $[R_{e,k}^{\lambda}]^{1/2}$ calculation, the state estimate computation still requires the $m\times m$ matrix inversion that is the inverse of $[R_{e,k}^{\lambda}]^{1/2}$; see line~\ref{mcc:sr1:p:X} of Algorithm~1. We stress that $[R_{e,k}^{\lambda}]^{1/2}$ is a triangular matrix and, hence, the inversion can be performed in efficient way (i.e., with the reduced computational cost) by solving the related linear equations through backward substitution. However, we can avoid the matrix inversion at all by the following trick. We augment the pre-array by a specially designed ``data'' column and apply the same orthogonal rotation to get the augmented post-array. Such algorithms are called  the {\it extended} implementations and they skip the matrix inversion for computing the state estimate.

Following~\cite{1995:Park}, consider the extended pre-array
\[
\begin{bmatrix}
R_k^{1/2} & \lambda_k^{1/2}H_{k}P_{k|k-1}^{1/2} & 0 \\
0 & F_kP_{k|k-1}^{1/2} & G_{k}Q_{k}^{1/2} \\
\hline
-\lambda_k^{1/2}y_k^{\top}R_k^{-{\top}/2} & \hat x_{k|k-1}^{\top}P_{k|k-1}^{-{\top}/2} & 0
\end{bmatrix}
\]
and lower-triangularize the first two (block) rows
\[\tilde {\mathbb A}_k {\mathbb Q}_k =
\begin{bmatrix}
[R_{e,k}^{\lambda}]^{1/2} & 0 & 0 \\
\bar K_{p,k}^{\lambda} & P_{k+1|k}^{1/2} & 0 \\
\hline
\alpha & \beta & \gamma
\end{bmatrix}.
\]

Similar to above presented derivation of Algorithm~1, we need to find $\alpha$, $\beta$ and $\gamma$. Taking into account properties of orthogonal matrices, from the inner product of the first and the last (block) rows of the pre-array, we get
\begin{align*}
 -&\lambda_k^{1/2} R_k^{1/2}R_k^{-1/2}y_k + \lambda_k^{1/2}H_{k}P_{k|k-1}^{1/2}P_{k|k-1}^{-1/2}\hat x_{k|k-1} \\
&  = [R_{e,k}^{\lambda}]^{1/2}\alpha^{\top} = -\lambda_k^{1/2}(y_k - H_k\hat x_{k|k-1}) = -\lambda_k^{1/2}e_k,
\end{align*}
i.e., $\alpha^{\top} = -\lambda_k^{1/2}[R_{e,k}^{\lambda}]^{-1/2}e_k$. We introduce notation for the normalized residuals
\[ \bar e_k = [R_{e,k}^{\lambda}]^{-1/2}e_k \mbox{ and } \bar e_k^{\lambda} =  \lambda_k^{1/2}\bar e_k = \lambda_k^{1/2}[R_{e,k}^{\lambda}]^{-1/2}e_k.
 \]
Hence, we conclude that $\alpha^{\top}=-\bar e_k^{\lambda}$. Next, equation~\eqref{MCCKF:X:Riccati1} of conventional one-step method is expressed in terms of available $[\bar K_{p,k}^{\lambda}]$ and the normalized residuals as follows:
\begin{align*}
\hat x_{k+1|k} & = F_k\hat x_{k|k-1} + \lambda_k^{1/2}[\bar K_{p,k}^{\lambda}]\bar e_k = F_k\hat x_{k|k-1} + [\bar K_{p,k}^{\lambda}]\bar e_k^{\lambda}
\end{align*}
where no matrix inversion is required for calculating the state, in contrast to the formula in line~\ref{mcc:sr1:p:X} of Algorithm~1.

It is worth noting here that the state vector can be calculated in alternative way. For that, we consider the inner product of the second and the last (block) rows of the pre-array in order to define $\beta$ value
\begin{align*}
\bar K_{p,k}^{\lambda}\alpha^{\top} + P_{k+1|k}^{1/2}\beta^{\top} & = -\bar K_{p,k}^{\lambda}\bar e_k^{\lambda} + P_{k+1|k}^{1/2}\beta^{\top}\\
&= F_kP_{k|k-1}^{1/2}P_{k|k-1}^{-1/2}\hat x_{k|k-1}  = F_k\hat x_{k|k-1},
\end{align*}
and, next, we conclude
\[\beta^{\top} = P_{k+1|k}^{-1/2}(F_k\hat x_{k|k-1}+\bar K_{p,k}^{\lambda}\bar e_k^{\lambda})=P_{k+1|k}^{-1/2}\hat x_{k+1|k}.\]

Thus, an alternative state computation way is to read-off values $[P_{k+1|k}^{-1/2}\hat x_{k+1|k}]$ and $[P_{k+1|k}^{1/2}]$ from the extended post-array, and simply multiply these blocks to obtain
\[\hat x_{k+1|k}= [P_{k+1|k}^{1/2}][P_{k+1|k}^{-1/2}\hat x_{k+1|k}].
 \]

In summary, the extended square-root condensed form implementation is presented in Algorithm~2.
\begin{codebox}
\Procname{{\bf Algorithm 2}. $\proc{eSR IMCC-KF}$ ({\it a priori, one-step, extended})}
\zi \textsc{Initialization:}($k=0$)
\li \>Apply Cholesky factorization: $\Pi_0 = \Pi_0^{1/2}\Pi_0^{{\top}/2}$;
\li \>Set initials: $P_{0|-1}^{1/2} = \Pi_0^{1/2}$, $P_{0|-1}^{-1/2}\hat x_{0|-1} =\Pi_0^{-1/2}\bar x_0$;
\zi \textsc{Filter Recursion}: ($k=\overline{0,K}$)
\li \>Compute $\lambda_k$ by equation~\eqref{eq:lambda};
\li \>Build and triangularize the first two (block) rows
\zi \>$\begin{array}{l}
\underbrace{\left[
\begin{smallmatrix}
R_k^{1/2} & \lambda_k^{1/2}H_{k}P_{k|k-1}^{1/2} & 0 \\
0 & F_kP_{k|k-1}^{1/2} & G_{k}Q_{k}^{1/2} \\
\hline
-\lambda_k^{1/2}y_k^{\top}R_k^{-{\top}/2^{\phantom{1/3}}} & [P_{k|k-1}^{-1/2}\hat x_{k|k-1}]^{\top} & 0
\end{smallmatrix}
\right]
}_{\mbox{\scriptsize Extended pre-array $\tilde {\mathbb A}_k$}} {\mathbb Q}_k \\
=
\underbrace{
\left[
\begin{smallmatrix}
[R_{e,k}^{\lambda}]^{1/2} & 0 & 0 \\
\bar K_{p,k}^{\lambda} & P_{k+1|k}^{1/2} & 0 \\
\hline
-[\bar e_k^{\lambda}]^{\top} & [P_{k+1|k}^{-1/2}\hat x_{k+1|k}]^{{\top}^{\phantom{{\top}/2}}} & (*)
\end{smallmatrix}
\right]
}_{\mbox{\scriptsize Extended post-array $\tilde {\mathbb R}_k$}};
\end{array}$
\li \>Read-off from post-array: $[P_{k+1|k}^{1/2}]$, $[P_{k+1|k}^{-1/2}\hat x_{k+1|k}]$;
\li \>Compute $\hat x_{k+1|k}= [P_{k+1|k}^{1/2}][P_{k+1|k}^{-1/2}\hat x_{k+1|k}]$.
\end{codebox}

Finally, we note that value $\gamma$ is of no interest for the filtering algorithm above. Hence, we use notation $(*)$ for the block that is not used in the filter. However, it is worth noting here that this quantity turns to be utilized in the square-root Rauch-Tung-Striebel formulas for the smoothed estimate $\hat x_{k|K}$, if one derives the smoother under the MCC methodology. In this case, the suggested Algorithm~2 yields a simple computation of this value.

\begin{table*}
\caption{The RMSE errors for IMCC-KF implementations in Example 1, $M = 100$ Monte Carlo simulations.} \label{tab:1}
\centering
{\small
\begin{tabular}{l|c|rrrr|r}
\hline
Estimator Implementation method & Covariance factorization & RMSE$_{x_1}$ & RMSE$_{x_2}$ & RMSE$_{x_3}$ & RMSE$_{x_4}$ & $\|$RMSE$\|_{2}$ \\
\hline
  & &\multicolumn{5}{c}{The errors are computed for {\it a priori} (predicted) estimates} \\
  \cline{3-7}
\texttt{IMCC-KF}: conventional eqs.~\eqref{MCCKF:P:Riccati1}--\eqref{MCCKF:X:Riccati1} &  --- & 6.7524	&  6.3572	& 38.7909	 & 38.5804	& 55.4905 \\
\texttt{SR IMCC-KF}: {\bf new} Algorithm~1  &  Cholesky & 6.7524	&  6.3572	& 38.7909	 & 38.5804	& 55.4905 \\
\texttt{eSR IMCC-KF}: {\bf new} Algorithm~2  & Cholesky & 6.7524	&  6.3572	& 38.7909	 & 38.5804	& 55.4905 \\	
\hline
  & & \multicolumn{5}{c}{The errors are computed for {\it a posteriori} (filtered) estimates} \\
  \cline{3-7}
\texttt{IMCC-KF}: eqs.~\eqref{mcckf:p:X}, \eqref{mcckf:p:P}, \eqref{mcckf:f:X}, \eqref{mcckf:f:K}--\eqref{mcckf:f:P} & --- & 6.6162	&  6.2099	& 38.7753	& 38.5790 &	 55.4455 \\
\texttt{SR IMCC-KF}: Algorithm~2 in~\cite{2017:CSL:Kulikova}	& Cholesky &   6.6162	&  6.2099 &	 38.7753 &	 38.5790	&  55.4455	 \\
\texttt{eSR IMCC-KF}: Algorithm~3 in~\cite{2017:CSL:Kulikova}	& Cholesky &  6.6162	&  6.2099 &	 38.7753 &	 38.5790	&  55.4455	 \\
\texttt{UD IMCC-KF}: Algorithm~2b in~\cite{2019:SP:Kulikova}	& modified Cholesky & 6.6162	&  6.2099 &	 38.7753 &	 38.5790	&  55.4455	 \\
\texttt{SVD IMCC-KF}: Algorithm~2c in~\cite{2019:SP:Kulikova}	& SVD &  6.6162	&  6.2099 &	 38.7753 &	 38.5790	&  55.4455	 \\
\hline
\end{tabular}
}
\end{table*}

\section{NUMERICAL EXPERIMENTS} \label{numerical:experiments}

The goal of this section is to substantiate the theoretical derivation and correctness   of the suggested condensed square-root Algorithms~1 and~2 on practical example.

\begin{exmp}[see~\cite{2016:Izanloo}] \label{ex:1} Consider a benchmark navigation problem where the vehicle dynamics is given as follows:
\begin{align*}
x_{k} & =
\begin{bmatrix}
1 & 0 & T & 0  \\
0 & 1 & 0 & T  \\
0 & 0 & 1 & 0 \\
0 & 0 & 0 & 1
\end{bmatrix}
x_{k-1} +
\begin{bmatrix}
w_{k-1}^1  \\
w_{k-1}^2  \\
w_{k-1}^3  \\
w_{k-1}^4
\end{bmatrix},
\end{align*}
where $T = 0.01$ is the sampling period, the first two state components are the north and east positions of a land vehicle,
and the last two components are the north and east velocities. A position-measuring device provides a noisy measurement of the north and
east positions as follows:
\begin{align*}
y_k & =
\begin{bmatrix}
1 & 0 & 0 & 0 \\
0 & 1 & 0 & 0
\end{bmatrix}
x_k +
\begin{bmatrix}
v_k^1 \\
v_k^2
\end{bmatrix}
\end{align*}
The MCC-KF methods are tested at the presence of impulsive noise, i.e., in the presence of outliers. More precisely, the process and measurement noises are generated as follows:
\begin{align*}
w_k  & \sim {\cal N}(0, Q)+\mbox{\tt Shot noise},  & Q & = 0.1 \; I_4\\
v_k  & \sim {\cal N}(0, R)+\mbox{\tt Shot noise},  & R & = 0.1 \; I_2
\end{align*}
where $\bar x_0 = [1,1,0,0]^T$ and $\Pi_0 = diag([4,4,3,3])$.
\end{exmp}

To simulate the impulsive noise (the shot noise), we follow the approach suggested in~\cite{2016:Izanloo}. The Matlab routine \verb"Shot_noise" recently published in~\cite[Appendix]{2019:AJC:Kulikova} can be used for generating the process and measurement noise in this set of experiments with 20\% of outliers distributed randomly in the time interval $[21,K-1]$ (where $K=300$) and with the randomly chosen  outliers' magnitude from $[0,5]$.

In our numerical simulations, we test all square-root IMCC-KF implementations existed nowadays~\cite{2017:CSL:Kulikova,2019:SP:Kulikova}. As mentioned in Introduction, they are  developed in two-step {\it a posteriori} form, only. The new one-step condensed Algorithms~1 and~2 are examined against {\it a priori}  one-step IMCC-KF conventional equations~\eqref{MCCKF:P:Riccati1}~--~\eqref{MCCKF:X:Riccati1}. All filtering algorithms under examination utilize the same initial filtering conditions, the same measurements and the same noise covariances. The experiment is repeated for $M=100$ Monte Carlo runs and the root mean square error (RMSE) is calculated for {\it a posteriori} and {\it a priori} estimates. The results of numerical experiments are summarized in Table~\ref{tab:1}.

We observe that all examined implementation methods produce the same RMSEs  both in {\it a priori} and {\it a posteriori} forms, respectively. In other words, they work with the same estimation accuracy. For  {\it a priori} one-step methods derived in this paper, this means that the results of practical example substantiate the correctness of their theoretical derivation presented in Section~\ref{main:result}. This also justifies an algebraic equivalence of new square-root condensed algorithms and the conventional way of calculations by equations~\eqref{MCCKF:P:Riccati1}~--~\eqref{MCCKF:X:Riccati1}.

\section{Conclusion} \label{conclusion}

This paper proposes one-step filtering strategy for the maximum correntropy criterion Kalman filters that corresponds to {\it a priori} filtering.
It has been shown that the resulted {\it a priori} recursion for the original MCC-KF does not allow for robust square-root implementations because of the possibly negative coefficient involved. In contrast, the improved MCC-KF estimator (IMCC-KF) easily permits stable
square-root implementations for its {\it a priori} recursion in the {\it condensed} form. Two Cholesky-based methods of this kind are derived in this paper. Finally, the square-root solution might be derived via $UDU^{\top}$ factorization as well. This problem is still open.

\bibliographystyle{IEEEtran}
\bibliography{IEEEabrv,BibTex_Library/books,%
              BibTex_Library/KFDistribution_Robust,%
              BibTex_Library/KF_Chandrasekhar,%
              BibTex_Library/KFMCC_Riccati,%
              BibTex_Library/KFMCC_Applications,%
              BibTex_Library/KFDiff_Chandrasekhar,%
              BibTex_Library/KF_Riccati,%
              BibTex_Library/CKF_Riccati,%
              BibTex_Library/UKF_Riccati,%
              BibTex_Library/EKF_all,%
              BibTex_Library/KF_Applications,%
              BibTex_Library/Lin_Algebra}

\end{document}